\documentclass[12pt,leqno]{article}
\usepackage[english]{babel}

\usepackage{amsfonts}
\usepackage{amssymb}
\usepackage{eucal}
\usepackage{amsxtra}
\usepackage{color}
\addtocounter{section}{-1} \numberwithin{equation}{section}
\newtheorem{theorem}{Theorem}[section]

\newtheorem{lemma}[theorem]{Lemma}

\newtheorem{corollary}[theorem]{Corollary}

\newcommand{\qed}{{$\hfill \Box$}}

\begin{document}

\title{\small \textbf{REAL HYPERSURFACES OF NON - FLAT COMPLEX SPACE FORMS IN TERMS OF THE JACOBI STRUCTURE OPERATOR}}
\bigskip
\author{\small Th. Theofanidis,  Ph. J. Xenos}
\date{}
\maketitle
\vspace{-25pt}
\begin{minipage}{400pt}
\begin{center}
Mathematics Division, School of Technology,\\
Aristotle University of Thessaloniki,\\
Thessaloniki, 54124, Greece.\\
email  :theotheo@gen.auth.gr, fxenos@gen.auth.gr.
\end{center}
\vspace{10pt}
\begin{abstract}\hspace{-18pt}The aim of the present paper is the study of
some classes of real hypersurfaces equipped with the condition $\phi l = l
\phi$, $(l = R( . , \xi)\xi).$
\end{abstract}
\
operator.
\end{minipage}

\section{Introduction.}

\hspace{15pt} An n - dimensional Kaehlerian manifold of constant holomorphic sectional
curvature c is called complex space form, which is denoted by $M_{n}(c)$.
The complete and simply connected complex space form is a projective space
$\mathbb{C}P^{n}$ if $c > 0$, a hyperbolic space $\mathbb{C}H^{n}$ if $c <
0$, or a Euclidean space $\mathbb{C}^{n}$ if $c = 0$. The induced almost
contact metric structure of a real hypersurface M of $M_{n}(c)$ will be
denoted by ($\phi, \xi, \eta, g$).

Real hypersurfaces in $\mathbb{C}P^{n}$ which are homogeneous, were
classified by R. Takagi (\cite{Takagi}). J. Berndt (\cite{Berndt})
classified  real hypersurfaces with principal structure vector fields in
$\mathbb{C}H^{n}$, which are divided into the model spaces $A_{0}$, $A_{1}$,
$A_{2}$ and $B$.

 Another class of real hypersurfaces were studied by
Okumura \cite{Okumura}, and Montiel and Romero \cite{Montiel Romero}, who
proved respectively the following theorems.
\begin{theorem}
Let M be a real hypersurface of  $\mathbb{C}P^{n}$, $n \geq 2$. If it
satisfies
$$ g((A \phi - \phi A)X, Y) = 0$$
for any vector fields X and Y, then M
is a tube of radius r over one of the following Kaehlerian submanifolds:\\
\
$(A_{1})$ a hyperplane $\mathbb{C}P^{n-1}$, where $0 < r <\frac{\pi}{2}$,\\
$(A_{2})$ a totally geodesic $\mathbb{C}P^{k}$($0 < k \leq n - 2$),where $0
< r <\frac{\pi}{2}$.
\end{theorem}

\begin{theorem}
Let M be a real hypersurface of  $\mathbb{C}H^{n}$, $n \geq 2$. If it
satisfies
$$ g((A \phi - \phi A)X, Y) = 0$$
for any vector fields X and Y, then M is locally congruent to
one of the following:\\
$(A_{0})$ a self - tube, that is, horosphere,\\
$(A_{1})$ a geodesic hypershere or a tube over a hyperplane $\mathbb{C}H^{n - 1}$,\\
$(A_{2})$ a tube over a totally geodesic $\mathbb{C}H^{k}$ ($1 \leq k \leq n
- 2$).
\end{theorem}

Real hypersurfaces of type $A_{1}$ and $A_{2}$ in $\mathbb{C}P^{n}$ and of
type $A_{0}$, $A_{1}$ and $A_{2}$ in $\mathbb{C}H^{n}$ are said to be
hypersurfaces of \emph{type  A} for simplicity.

A Jacobi field along geodesics of a given Riemannian manifold (M, g) plays
an important role in the study of differential geometry. It satisfies a well
known differential equation which inspires Jacobi operators. For any vector
field $X$, the Jacobi operator is defined by $R_{X}$: $R_{X}(Y) = R(Y,
X)X$, where $R$ denotes the curvature tensor and $Y$ is a vector field on
M. $R_{X}$ is a self - adjoint endomorphism in the tangent space of M, and
is related to the Jacobi differential equation, which is given by
$\nabla_{\acute{\gamma}}(\nabla_{\acute{\gamma}}Y) + R(Y,
\acute{\gamma})\acute{\gamma} = 0$ along a geodesic $\gamma$ on $M$, where
$\acute{\gamma}$ denotes the velocity vector along $\gamma$ on $M$.

In a real hypersurface $M$ of a complex space form $M_{n}(c)$, $c \neq 0$,
 the Jacobi operator on $M$ with respect to the structure vector field $\xi$,
 is called the  Jacobi structure operator and is denoted by
$l X = R_{\xi}(X) = R(X, \xi)\xi$.

Many authors have studied real hypersurfaces from many points of view.
Certain authors have studied real hypersurfaces under the condition $\phi l
= l \phi$, equipped with one or two additional conditions. U-Hang Ki, An -Aye Lee and Seong-Baek Lee (\cite{Ki Lee Lee}) classified real hypersurfaces in complex space forms satisfying i) $\phi l= l \phi$ and $A^{2}\xi = \theta A \xi + \tau \xi$ ($\theta$ is a function, $\tau$ is constant) ii) $\phi l= l \phi$  and $Q\xi = \sigma\xi$ (where $Q$ is the Ricci operator, $\sigma$ is constant). U-Hang Ki (\cite{Ki}) classified real hypersurfaces in complex hyperbolic space satisfying $\phi l= l \phi$ and $lQ = Ql$. U-Hang Ki with Soo Jin-Kim and Seong-Baek Lee (\cite{Ki Kim Lee}), classified real hypersurfaces in complex space forms satisfying $\phi l= l \phi$, $lQ = Ql$, and additional conditions on the mean curvature. U-Hang Ki, S. Nagai and R. Takagi(\cite{Ki Nagai Takagi}) studied real hypersurfaces in complex space forms satisfying $\phi l= l \phi$ and $l Q= l Q$.

Other authors have studied real hypersurfaces under the conditions  $\nabla_{X}l = 0$ ($X \in TM$) or  $\nabla_{\xi}l = 0$ (\cite{Ortega Perez Santos}, \cite{Ki Perez Santos}, \cite{Cho Ki1}).

In the present paper, we consider a weaker condition $\nabla_{\xi}l = \mu\xi$, where $\mu$ is a function
of class $C^{1}$ on M, and classify these hypersurfaces satisfying $\phi l = l \phi$. Namely we prove:

\begin{theorem}
Let M be a real hypersurface of a complex space form $M_{n}(c)$, ($n > 2$) $(c \neq
0)$, satisfying
$\phi l = l\phi$.  If  $\nabla_{\xi}l = \mu \xi$ on ker($\eta$) or on  span$\{\xi\}$, then M is  a Hopf hypersurface. Furthermore, if $\eta(A\xi) \neq 0$, then M locally congruent to a model space of type A.
\end{theorem}

J. T. Cho and U - H Ki in \cite{Cho Ki} classified real hypersurfaces M of a projective space
satisfying $\phi l = l\phi$ and $lA = Al$ on M. In the
 present paper we generalize this result, studying the real
 hypersurfaces of any complex space form satisfying $\phi l = l\phi$ and $lA = Al$ on the distribution on M (ker($\eta$)) given by all vectors orthogonal to the Reeb flow $\xi$, or on span$\{\xi\}$. We prove:
\begin{theorem}
Let M be a real hypersurface of a complex space form $M_{n}(c)$, ($n > 2$) $(c \neq
0)$, satisfying
$\phi l = l\phi$.  If  $l A = A l$ on ker($\eta$) or on span$\{\xi\}$, then M is  a Hopf hypersurface. Furthermore, if $\eta(A\xi) \neq 0$, then M locally congruent to a model space of type A.
\end{theorem}

For the case of $\mathbb{C}P^{n}$ in order
to determine real hypersurface of type A, the
technical assumption $\eta(A\xi) \neq 0$ is needed.
Actually, there is a non-homogeneous tube with $A\xi = 0$ (of radius $\frac{\pi}{4}$)
 over a certain Kaehler submanifold in $\mathbb{C}P^{n}$,
 when its focal map has constant rank on $M$ (\cite{Cecil Ryan}). For
  Hopf hypersurfaces in $\mathbb{C}H^{n}$, $(n > 2)$ it is known that the
  associated principal curvature of $\xi$ never vanishes (\cite{Berndt}).
   However, in $\mathbb{C}H^{2}$ there exists a Hopf hypersurface
   with $A\xi = 0$ (\cite{Ivey Ryan}).
\section{Preliminaries.}

\hspace{15pt} Let $M_{n}$ be a Kaehlerian manifold of real
dimension 2n, equipped with an almost complex structure J and a Hermitian
metric tensor G. Then for any vector fields X and Y on $M_{n}(c)$, the
following
relations hold:\\
\begin{center}
$J^{2}X = -X$, \hspace{20pt}$G(JX, JY) = G(X, Y)$,\hspace{20pt}
$\widetilde{\nabla}J = 0$\end{center} where $\widetilde{\nabla}$ denotes the
Riemannian connection of G of
$M_{n}$.\\

Now, let $M_{2n-1}$ be a real (2n-1)-dimensional hypersurface of $M_{n}(c)$,
and denote by N a unit normal vector field on a neighborhood of a point in
$M_{2n-1}$ (from now on we shall write \emph{M} instead of $M_{2n-1}$). For
any vector field  X tangent to M we have $JX = \phi X + \eta(X)N$, where
$\phi X$ is the tangent component of $JX$, $\eta(X)N$ is the normal
component, and
\begin{center}
$\xi = - JN$, \hspace{30pt} $\eta(X) = g(X, \xi)$, \hspace{30pt}$g =
G|_{M}$.
\end{center}

By properties of the almost complex structure J, and the definitions of
$\eta$ and g, the following relations hold (\cite{Blair}):

\begin{center}
\begin{flushleft}
(1.1)\hspace{70pt}$\phi^{2} = - I + \eta\otimes\xi, \hspace{20pt}
\eta\circ\phi = 0 ,
    \hspace{20pt} \phi\xi = 0, \hspace{10pt}
 \eta(\xi) = 1 $\newline
\ \\ (1.2)\hspace{20pt}$g(\phi X, \phi Y) = g(X, Y) - \eta(X)\eta(Y),
\hspace{20pt} g(X, \phi Y) = - g(\phi X, Y)$.\end{flushleft}

\end{center}
The above relations define an \emph{almost contact metric structure} on M
which is denoted by $(\phi, \xi, g, \eta)$. When an almost contact metric
structure is defined on M, we can define a local orthonormal basis $\{V_{1},
V_{2}, . . . V_{n - 1}, \phi V_{1}, \phi V_{2}, . . . \phi V_{n - 1},
\xi\}$, called a $\phi - basis$. Furthermore, let A be the shape operator in
the direction of N, and denote by $\nabla$ the Riemannian connection of g on
M. Then, A is symmetric and the following equations are
satisfied:\\
\begin{center}
\begin{flushleft}
(1.3) \hspace{40pt} $\nabla_{X}\xi = \phi AX  ,\hspace{40pt}
(\nabla_{X}\phi)Y = \eta(Y)AX -
    g(AX, Y)\xi$.\end{flushleft}
\end{center}

As the ambient space $M_{n}(c)$ is of constant holomorphic sectional
curvature c, the equations of Gauss and Godazzi are respectively given by:

\begin{center}
(1.4) \hspace{20pt}$R(X, Y)Z = \frac{c}{4}[g(Y, Z) X - g(X, Z)Y + g(\phi Y,
Z)\phi X - g(\phi
X, Z)\phi Y$\\
\ \\
\hspace{20pt}$ - 2g(\phi X, Y)\phi Z] + g(AY, Z)AX - g(AX,
Z)AY$,\end{center} \
\begin{center}
(1.5)\hspace{60pt}$ (\nabla_{X}A)Y - (\nabla_{Y}A)X =
\frac{c}{4}[\eta(X)\phi Y - \eta(Y)\phi X - 2g(\phi X, Y)\xi]$.
\end{center}

The tangent space $T_{p}M$, for every point $p\in M$, is decomposed as
following:
\begin{center}
$T_{p}M = ker(\eta)^{\bot} \oplus ker(\eta)$
\end{center}
where $ker(\eta)^{\bot} = span\{\xi\}$ and $ker(\eta)$ is defined as following:\\
\begin{center}
$ker(\eta) = \large{\{}
X \in T_{p}M: \eta(X) = 0\large{\}}$ \\
\end{center}
\
\\
Based on the above decomposition, by virtue of (1.3), we decompose the
vector field $A\xi$ in the following way:
\begin{flushleft}
(1.6) \hspace{100pt}$A\xi = \alpha \xi + \beta U$\end{flushleft}

\begin{flushleft}
where $\beta = |\phi\nabla_{\xi}\xi|$ and $U = -\frac{1}{\beta}
\phi\nabla_{\xi}\xi \in ker(\eta)$, provided that $\beta \neq 0$.
\end{flushleft}

If the vector field $A\xi$ is expressed as $A\xi = \alpha \xi$, then $\xi$
is called a principal vector field.

Finally differentiation will be denoted by ( ). All manifolds and vector fields of this paper
are assumed to be connected and of class $C^{\infty}$.

\section{Auxiliary relations}

\hspace{15pt}In the study of real hypersurfaces of a complex space form $M_{n}(c)$, $c
\neq 0$, it is a crucial condition that the structure vector field $\xi$ is
principal.  The purpose of this paragraph is
to prove this condition.

Let V be the open subset of points p of M, where $\alpha\neq0$ in a
neighborhood of p and $V_{0}$ be the open subset of points p of M such that
$\alpha = 0$ in a neighborhood of p. Since $\alpha$ is a smooth function on
M, then $V\cup V_{0}$ is an open and dense subset of M.

\begin{lemma}
Let M be a real hypersurface of a complex space form $M_{n}(c)$ $(c \neq
0)$, satisfying $\phi l = l \phi$ on $ker(\eta).$ Then, $\beta = 0$ on
$V_{0}$.
\end{lemma}
\vspace{-15pt}
Proof. From (1.6) we have $A\xi = \beta U$ on $V_{0}$. Then (1.4) for $X = U$ and
$Y = Z= \xi$ yields
$$lU = \frac{c}{4}U + g(A\xi, \xi)AU - g(AU, \xi)A\xi = \frac{c}{4}U
- g(U, A\xi)A\xi = (\frac{c}{4} - \beta^{2})U \Rightarrow$$
$$\phi l U = (\frac{c}{4} - \beta^{2})\phi U. $$
In the same way, from (1.4) for $X = \phi U$, $Y = Z = \xi$ we obtain
$$l \phi U = \frac{c}{4} \phi U. $$
The last two equations yield $\beta = 0.$ \qed\\
\emph{\textbf{REMARK 1}}\\
We have proved that on $V_{0}$, $A\xi = 0 \xi$ i.e. $\xi$ is a principal
vector field on $V_{0}$. Now we define on V the set $V'$ of points p where
$\beta \neq 0$ in a neighborhood of p and the set $V''$ of points p where
$\beta = 0$ in a neighborhood of p. Obviously $\xi$ is principal on $V''$.
In what follows we study the open subset $V'$ of M and define the following classes:\\
A =  hypersurfaces satisfying $\phi l = l \phi$ and $l A = A l$ on ker{($\eta$)},\\
B =  hypersurfaces satisfying $\phi l = l \phi$ and $l A = A l$ on span\{$\xi$\},\\
C =  hypersurfaces satisfying $\phi l = l \phi$ and $\nabla_{\xi} l = \mu\xi$ on ker{($\eta$)},\\
D =  hypersurfaces satisfying $\phi l = l \phi$ and $\nabla_{\xi} l = \mu\xi$ on on span\{$\xi$\}.\\

\begin{lemma}
Let M be a real hypersurface of a complex space form $M_{n}(c)$
($c\neq 0$), satisfying $\phi l = l \phi$ on $ker(\eta)$ . Then the following relations hold on the set $V'$ of classes A, B , C, D.\\
\end{lemma}
\begin{equation}
AU = (\frac{\beta^{2}}{\alpha} - \frac{c}{4\alpha}
)U + \beta \xi, \hspace{50pt}A \phi U = - \frac{c}{4\alpha}\phi U.
\end{equation}

\begin{equation}
\nabla_{\xi}\xi = \beta \phi U , \hspace{20pt}  \nabla_{U}\xi =
(\frac{\beta^{2}}{\alpha} - \frac{c}{4\alpha}) \phi U ,\hspace{20pt} \nabla_{\phi U}\xi =
 \frac{c}{4\alpha} U.
\end{equation}

\begin{equation}
\nabla_{\xi}U = W_{1}  , \hspace{20pt} \nabla_{U}U = W_{2}, \hspace{20pt} \nabla_{\phi U}U = W_{3} - \frac{c}{4\alpha}\xi.
\end{equation}

\begin{equation}
\nabla_{\xi}\phi U = \phi W_{1} - \beta \xi, \hspace{10pt}\nabla_{U}\phi
U = \phi W_{2} + (\frac{c}{4\alpha} -
\frac{\beta^{2}}{\alpha})\xi, \nabla_{\phi U}\phi U = \phi W_{3}.
\end{equation}
where $W_{1}$, $W_{2}$, $W_{3}$ are vector fields on ker$(\eta)$ satisfying $W_{1}, W_{2}, W_{3}\perp U$.\\
\
\\
Proof. From (1.4) we get \
\\
\begin{equation}
    lX = \frac{c}{4}[X - \eta(X)\xi] + \alpha AX - g(AX, \xi)A\xi
\end{equation}
which, for X = U yields
\begin{equation}
    lU = \frac{c}{4}U + \alpha AU - \beta A\xi.
\end{equation}
The scalar products of (2.6) with U and $\phi U$ yield respectively
\begin{equation}
    g(AU, U) = \frac{\gamma}{\alpha} - \frac{c}{4\alpha} + \frac{\beta^{2}}{\alpha},
\end{equation}
\begin{equation}
g(AU, \phi U) = \frac{1}{a}g(lU, \phi U).
\end{equation}
where $\gamma = g(lU, U) = g(\phi l U, \phi U) = g(l\phi U, \phi U).$\\
The second relation
of (1.2) for $X = U$, $Y = lU$, the condition $\phi l =l \phi$ and the
symmetry of the operator $l$ imply:
\begin{center}
$g(lU, \phi U) = 0 $.
\end{center}
The above equation and (2.8) imply
\begin{equation}
g(AU, \phi U) = 0.
\end{equation}
The symmetry of $A$ and (1.6) imply
\begin{equation}
g(AU, \xi) = \beta.
\end{equation}
From relations (2.7), (2.9) and (2.10), we obtain

\begin{equation}AU = (\frac{\gamma}{\alpha}-\frac{c}{4\alpha} +
\frac{\beta^{2}}{\alpha})U + \beta \xi +
\lambda W
\end{equation}
 where $W\in span\{U, \phi U, \xi \}^{\bot}$ and $\lambda = g (AU, W)$.
Combining (2.11) with (2.6) we obtain $l U = \gamma U + \lambda \alpha W$.
Acting on this relation with the tensor field $\phi$ and by virtue
of $\phi l = l\phi$ we take $l \phi U = \gamma \phi U + \lambda \alpha \phi W$.
On the other hand by virtue of (2.5) we have $l\phi U = \frac{c}{4}\phi U + \alpha A\phi U$.
 From the last two relations we obtain  $A\phi U = (\frac{\gamma}{\alpha}-\frac{c}{4\alpha})\phi U + \lambda \phi W.$\\
\underline{On class A}\\
Since $lA = Al$ holds on ker($\eta$) we have $lAW = AlW$.
This relation because of (2.5) and (2.11) implies $\lambda\beta A\xi = 0$ and
so $\lambda = 0$. Since $\lambda = 0$, equations $lAU = AlU$, (2.6) and (2.11)
 yield $\gamma = 0$, therefore we have the first of (2.1). Moreover from
 (2.5) we have $l\phi U = \frac{c}{4}\phi U + \alpha A\phi U$ which is written as
 $\phi l U =  \frac{c}{4}\phi U + \alpha A\phi U$ ($\phi l = l\phi$). From $\phi l U =  \frac{c}{4}\phi U + \alpha A\phi U$ and $\gamma = \lambda = 0$ we obtain the second of (2.1). Using (1.3) for $X \in \{\xi, U, \phi U\}$  and by virtue of (2.1) we obtain (2.2).
It is well known that:
\begin{equation}
X g(Y, Z) = g(\nabla_{X}Y, Z) +   g(Y, \nabla_{X}Z)
\end{equation}
Let us set $\nabla_{\xi}U = W_{1}$ and $\nabla_{U}U = W_{2}$. If we use (2.2) and (2.12), it is easy to verify that $g(\nabla_{\xi}U, U) = 0 = \eta(\nabla_{\xi}U)$ and $g(\nabla_{U}U, U) = 0 = \eta(\nabla_{U}U) $ which means $W_{1}\bot \{\xi,U\}$ and $W_{2}\bot \{\xi,U\}$.

On the other hand using (2.12) and the third of (2.2) we
 find $\eta(\nabla_{\phi U}U) = -\frac{c}{4\alpha}$ and $g(\nabla_{\phi U}U, U) = 0$
  which means that $\nabla_{\phi U}U$ is decomposed as
   $\nabla_{\phi U}U = W_{3} -\frac{c}{4\alpha}\xi$, $W_{3}\bot \{U, \xi \}.$
Now, by virtue of (1.3) and (2.3) for $X = \xi, Y = U$ and $X = Y =
U$ and
 $X = \phi U, Y = U$, we get (2.4).
\\
\
\\
\underline{On class B}\\
We analyze equation $l A \xi= Al\xi$ by virtue of (1.6), (2.6) and (2.11) and we have $\gamma U + \lambda\alpha W = 0$. Since $W\perp U$ we have $\gamma = \lambda = 0$. The rest of the proof is similar to the one in class A.
\
\\
\underline{On class C}\\
We have $(\nabla_{\xi}l)U = \mu\xi$. The scalar product of
$(\nabla_{\xi}l)U = \mu\xi$ with $\xi$, the symmetry of $l$, $g(lU,
\phi U) = 0$ and (2.12) yield $\mu = 0$. In addition we have
$(\nabla_{\xi}l)\phi U = \mu\xi = 0$. So, the scalar product of
$(\nabla_{\xi}l)\phi U = 0$ with $\xi$,
the symmetry of $l$ and (2.12) yield $g(l \phi U, \phi U) = \gamma = 0.$\\
 Finally $(\nabla_{\xi}l)\phi W = \mu\xi = 0$. So, the scalar product of $(\nabla_{\xi}l)\phi W = 0$ with $\xi$, the symmetry of $l$  and (2.12) yield $g(l\phi U, \phi W) = 0$, which, by virtue of (2.5), the second of (2.1) and $\gamma = 0$, yield $\lambda = 0$. The rest of the proof is similar to the one in class A.\\
 \underline{On class D}\\
We analyze $(\nabla_{\xi}l)\xi = \mu\xi$ and obtain $\beta l\phi U = \mu\xi$.
But the vector fields $l\phi U$ and $\xi$ are linear independent, so
 $l\phi U = \mu = 0$. We analyze $l\phi U = 0$ using (2.5) and  $A\phi U =
 (\frac{\gamma}{\alpha}-\frac{c}{4\alpha})\phi U + \lambda \phi W$,
 and we have $\gamma \phi U + \lambda\alpha\phi W = 0$. This relation and the linear independency
 of the vector fields $\phi U$ and $\phi W$ yield $\gamma = \lambda = \mu = 0$.
  The rest of the proof is similar to the one in class A.
\qed

\begin{lemma}
Let M be a real hypersurface of a complex space form $M_{n}(c)$ $(c \neq
0)$, of class A, B ,C, or D. Then on $V'$ we have
$g(\nabla_{\xi}U, \phi U) = -4\alpha$ and $g(\nabla_{U}U, \phi U) = -4\beta + \frac{c}{4\alpha\beta}(\frac{c}{4\alpha} - \frac{\beta^{2}}{\alpha})$.\end{lemma}
 Proof. Putting $X = U$, $Y = \xi$ in (1.5), we obtain
$$(\nabla_{U}A)\xi - (\nabla_{\xi}A)U = - \frac{c}{4}\phi U.$$
Combining the last equation with (1.6), and Lemma 2.2
 it follows :
$$(U\alpha) \xi + (U\beta)U + \beta W_{2} + (-\frac{c}{4\alpha} + \frac{\beta^{2}}{\alpha})\frac{c}{4\alpha}\phi U$$ $$- \xi(-\frac{c}{4\alpha}+ \frac{\beta^{2}}{\alpha})U - (-\frac{c}{4\alpha} + \frac{\beta^{2}}{\alpha})W_{1} - (\xi\beta)\xi + AW_{1} = 0.$$

Taking the scalar products of the last relation with $\xi$ and $U$ respectively, we obtain
\begin{equation}
(U\alpha) = (\xi\beta)
\end{equation}and
\begin{equation}
(U\beta) = (\xi(\frac{\beta^{2}}{\alpha} - \frac{c}{4\alpha})).
\end{equation}
Combining the last three equations we have
\begin{equation}
AW_{1} =  \frac{c}{4\alpha}(\frac{c}{4\alpha} - \frac{\beta^{2}}{\alpha})\phi U + (\frac{\beta^{2}}{\alpha} - \frac{c}{4\alpha})W_{1} - \beta W_{2}.
\end{equation}
The scalar product of (2.15) with $\phi W_{1}$ yields: $$\beta
g(\phi W_{1}, W_{2}) = - g(AW_{1}, \phi W_{1}).$$ But from (2.5) we
have $$g(l\phi W_{1}, W_{1}) = g(\phi W_{1}, l W_{1}) = \alpha g(A
W_{1}, \phi W_{1}).$$ Moreover $ g(l\phi W_{1}, W_{1}) = g(\phi
W_{1}, l W_{1}) = - g( W_{1}, \phi l W_{1}) =
- g( W_{1},l \phi W_{1})$ which means that $$g(l\phi W_{1}, W_{1}) = 0.$$
 The above relations lead to $g(\phi W_{1}, W_{2}) = 0$ which, by virtue of (2.15) implies $g(A W_{1}, \phi W_{2}) = 0$.\\
In what follows we define the following functions: $$\kappa_{1} = g(W_{1}, \phi U)
\hspace{20pt} \kappa_{2} = g(W_{2}, \phi U), \hspace{20pt}\kappa_{3} = g(W_{3}, \phi U).$$

Putting $X = \phi U$, $Y = \xi$ in (1.5), we obtain
\begin{equation}
 A\phi W_{1} =  [ \frac{3\beta c}{4\alpha} + \alpha\beta -(\phi U \alpha)]\xi\end{equation} $$-
 [(\phi U \beta) + \frac{c}{4\alpha}(\frac{c}{4\alpha} - \frac{\beta^{2}}{\alpha}) -
  \beta^{2}]U  + \frac{c}{4\alpha^{2}}(\xi \alpha)\phi U - \frac{c}{4\alpha}\phi W_{1} - \beta
  W_{3}.$$

The scalar product of (2.16) with $\xi$ implies
\begin{equation}
(\phi U \alpha) = \frac{3\beta c}{4\alpha} + \alpha\beta +
\kappa_{1}\beta.
\end{equation}
From the scalar product of (2.16) with $U$ we get
$$
g(A\phi W_{1}, U) = -(\phi U \beta) -
\frac{c}{4\alpha}(\frac{c}{4\alpha} - \frac{\beta^{2}}{\alpha}) +
\beta^{2} - \frac{c}{4\alpha}g(\phi W_{1}, U) \Rightarrow$$
$$
g(\phi W_{1}, A U) = -(\phi U \beta) -
\frac{c}{4\alpha}(\frac{c}{4\alpha} - \frac{\beta^{2}}{\alpha}) +
\beta^{2} + \frac{c}{4\alpha}g(W_{1}, \phi U) \Rightarrow,$$ which,
eventually (with the aid of Lemma 2.2 and the definition of
$\kappa_{1}$) yields
\begin{equation}
(\phi U \beta) =  \frac{c}{4\alpha}(\frac{\beta^{2}}{\alpha} - \frac{c}{4\alpha}) + \beta^{2} + \kappa_{1}\frac{\beta^{2}}{\alpha}.
\end{equation}

The condition $\phi l W_{1} = l\phi W_{1}$ because of (2.5), (2.15) and (2.16) implies
$$
\beta^{2}\phi W_{1} - \alpha\beta\phi W_{2} + \alpha[(\phi U \alpha) - \frac{3\beta c}{4\alpha} -
\alpha\beta]\xi + \alpha[(\phi U \beta) - \beta^{2}]U + \alpha\beta W_{3}$$ $$= \kappa_{1}\beta A\xi + \frac{c}{4\alpha}(\xi\alpha)\phi U.$$
\
\\
Taking the scalar product of the last relation with $U$ we have
$$
-2\kappa_{1}\beta^{2} + \alpha\beta\kappa_{2} + \alpha(\phi U \beta) - \alpha\beta^{2} = 0.$$
If in the above relation we replace the term $\kappa_{1}$ using (2.17) we obtain
\begin{equation}
-2\beta(\phi U \alpha) + \frac{3\beta^{2}c}{2\alpha} + \alpha\beta^{2} + \alpha\beta\kappa_{2} + \alpha(\phi U \beta) = 0.
\end{equation}

The relation $(\nabla_{U}A)\phi U - (\nabla_{\phi U}A) U = -\frac{c}{2}\xi$, using Lemma 2.2 implies
\begin{equation}
\frac{c}{4\alpha^{2}}(U \alpha) \phi U + [\frac{c}{2\alpha}( \frac{\beta^{2}}{\alpha} -\frac{c}{4\alpha}) + \beta^{2} - (\phi U \beta)]\xi +\end{equation}

$$ [-\frac{3\beta c}{4\alpha} + \frac{\beta^{3}}{\alpha} + (\phi U(\frac{c}{4\alpha} - \frac{\beta^{2}}{\alpha}))]U
- \frac{c}{4\alpha} \phi W_{2} - A\phi W_{2}$$ $$+ AW_{3} + (\frac{c}{4\alpha} -\frac{\beta^{2}}{\alpha})W_{3} = 0.$$
The scalar product of the above relation with U yields

 $$\frac{\kappa_{2}\beta^{2}}{\alpha} - \frac{3\beta c}{4\alpha} + \frac{\beta^{3}}{\alpha}
  + \phi U(\frac{c}{4\alpha} - \frac{\beta^{2}}{\alpha}) = 0.$$
Expanding the last relation and by virtue of (2.19) we get
$$
(- \frac{3\beta^{2}}{\alpha^{2}} + \frac{c}{4\alpha^{2}})(\phi U \alpha) + \frac{3\beta}{\alpha}(\phi U\beta) + \frac{3\beta^{3}c}{2\alpha^{3}} + \frac{3\beta c}{4\alpha} = 0.
$$

Combining the last equation with (2.17) and (2.18) we obtain $\kappa_{1} = -4\alpha.$ The scalar product of (2.15) with $\phi U$  because of $\kappa_{1} = -4\alpha$, yields $\kappa_{2} = -4\beta + \frac{c}{4\alpha\beta}(\frac{c}{4\alpha} - \frac{\beta^{2}}{\alpha}).$\qed

\begin{lemma}
Let M be a real hypersurface of a complex space form $M_{n}(c)$ $(c \neq
0)$, of class A, B ,C, or D. Then the structure vector field $\xi$ is principal
on M.
\end{lemma}
Proof.  The scalar products of (2.16) and (2.20) with $\phi U$,
yield $(\xi\alpha) = \frac{4\alpha^{2}\beta}{c}\kappa_{3}$ and
$(U\alpha) = \frac{4\alpha\beta^{2}}{c}\kappa_{3} $. Combining the
last two relations with (2.13) and (2.14) we have
\begin{equation}(\xi\alpha) = \frac{4\alpha^{2}\beta}{c}\kappa_{3},\hspace{10pt} (U\alpha) = (\xi\beta) =
\frac{4\alpha\beta^{2}}{c}\kappa_{3}, \hspace{10pt} (U\beta) = (\beta + \frac{4\beta^{3}}{c})\kappa_{3}
\end{equation}

Using (1.5) for $X = \phi W_{2}$, $Y = \xi$ we have

$$\nabla_{\phi W_{2}}A\xi - A\nabla_{\phi W_{2}}\xi - \nabla_{\xi}A\phi W_{2} + A\nabla_{\xi}\phi W_{2} = \frac{c}{4} W_{2},$$
which, from (1.6) is further decomposed as
$$(\phi W_{2}\alpha)\xi + \alpha \phi A\phi W_{2} + (\phi W_{2}\beta) U  +  \beta\nabla_{\phi W_{2}}U -$$ $$A\phi A\phi W_{2} - \nabla_{\xi}A\phi W_{2} + A\nabla_{\xi}\phi W_{2} = \frac{c}{4} W_{2}.$$
Taking the scalar product with $\xi$ and by using (1.6), (2.12), (2.21), Lemmas 2.2, 2.3 and $W_{1}\bot\phi W_{2}$ we obtain
\begin{equation}
(\phi W_{2}\alpha) = \kappa_{3}(\frac{16 \alpha\beta^{3}}{c} +
\beta(\frac{\beta^{2}}{\alpha} - \frac{c}{4\alpha}) ).
\end{equation}

On the other hand from (1.5) we get
$$\nabla_{ W_{3}}A\xi - A\nabla_{W_{3}}\xi - \nabla_{\xi}A W_{3} + A\nabla_{\xi} W_{3} = -\frac{c}{4}\phi W_{3}$$
which, by virtue of (1.6) is further decomposed as
$$(W_{3}\alpha)\xi + \alpha\phi A W_{3} + (W_{3}\beta)U + \beta\nabla_{W_{3}}U - A\nabla_{W_{3}}\xi
 - \nabla_{\xi}A W_{3} + A\nabla_{\xi} W_{3} = -\frac{c}{4}\phi W_{3}.$$
Taking the scalar product of the last equation with $\xi$ and by
making use of Lemma 2.2, (2.12) and (2.21)we obtain
\begin{equation}
(W_{3}\alpha) = 3\beta(\frac{c}{4\alpha} - \alpha)\kappa_{3}.
\end{equation}

In a similar way equation (1.5) yields $(\nabla_{\phi W_{1}}A)U -
(\nabla_{U}A)\phi W_{1} = 0$,which by virtue of Lemma 2.2 is further
analyzed as $$(\phi W_{1}(\frac{\beta^{2}}{\alpha} - \frac{c}{4\alpha}))U +
(\frac{\beta^{2}}{\alpha} - \frac{c}{4\alpha})\nabla_{\phi W_{1}}U +$$ $$
(\phi W_{1}\beta)\xi + \beta\phi A\phi W_{1} - A\nabla_{\phi W_{1}}U
- \nabla_{U}A\phi W_{1} + A\nabla_{U}\phi W_{1} = 0.$$ The scalar
product of the above equation with $\xi $ and using $g(\phi W_{1}, W_{2}) = 0$, (2.21) and Lemmas 2.2,  2.3, leads to
\begin{equation}
(\phi W_{1}\beta) = 4\alpha\kappa_{3}(\beta + \frac{4\beta^{3}}{c}).
\end{equation}

Now the calculation of Lie bracket $[\phi U, \xi]\beta$, by virtue
of (2.18), Lemma 2.3 and (2.21), results to
$$[\phi U, \xi]\beta = \phi U(\xi\beta) + \kappa_{3}[-\frac{\beta c}{2\alpha} + \frac{24\alpha\beta^{3}}{c}]. $$

On the other hand from Lemma 2.2, (2.21) and (2.24) we obtain
$$[\phi U, \xi]\beta = (\nabla_{\phi U}\xi - \nabla_{\xi}\phi U)\beta = \beta\kappa_{3}[\frac{c}{4\alpha} + \frac{\beta^{2}}{\alpha} - 4\alpha
 - \frac{12\alpha\beta^{2}}{c}].$$

Equalizing the above two relations we get
\begin{equation}
\phi U(\xi\beta) = \beta\kappa_{3}[\frac{3c}{4\alpha} +
\frac{\beta^{2}}{\alpha} - 4\alpha
 - \frac{36\alpha\beta^{2}}{c}].
\end{equation}

In a similar way, combining (2.18), (2.21), (2.22), (2.23) and Lemmas 2.2, 2.3, the Lie
bracket $[\phi U, U]\alpha$
yields
$$[\phi U, U]\alpha = \phi U(U\alpha) + 3\beta\kappa_{3}[\alpha
 + \frac{8\alpha\beta^{2}}{c} - \frac{c}{4\alpha}].$$
$$
[\phi U, U]\alpha = (\nabla_{\phi U} U - \nabla_{U}\phi U)\alpha = \beta \kappa_{3}[\frac{c}{\alpha} - 5\alpha
 - \frac{12\alpha\beta^{2}}{c} - \frac{\beta^{2}}{\alpha}].
$$
From the above equations we obtain
\begin{equation}
\phi U( U\alpha) = \beta\kappa_{3}[\frac{7 c}{4\alpha} -
8\alpha
 - \frac{36\alpha\beta^{2}}{c} - \frac{\beta^{2}}{\alpha}].
\end{equation}

Because of (2.13) from (2.25) and (2.26)we obtain
$$\frac{\beta}{\alpha}[c - 4\alpha^{2} - 2
\beta^{2}]\kappa_{3} = 0.$$

Let us assume there is a point $p\in V'$
such that $\kappa_{3}\neq 0$  in a neighborhood around $p$. Then we
have $c = 4\alpha^{2} + 2 \beta^{2}.$ Differentiating the last
equation along $\xi$ and by virtue of (2.21) and $\kappa_{3} \neq 0$
we take $2\alpha^{2} + \beta^{2} = 0$ which is a contradiction. So
$\kappa_{3} = 0 \Rightarrow (U\alpha) = (\xi\alpha) = 0 \Rightarrow
[U, \xi]\alpha = 0.$ But the last equation, because of Lemma 2.2
yields
\begin{equation}
(\frac{\beta^{2}}{\alpha} - \frac{c}{4\alpha}) (\phi U \alpha) - (W_{1}\alpha)= 0.
\end{equation}

On the other hand from (1.5) for $X = W_{1}$, $Y = \xi$, taking the scalar product with $\xi$, using the Lemmas 2.2, 2.3 we have $(W_{1}\alpha) = \beta|W_{1}|^{2} - \beta(4\alpha^{2} + 3 c)$. The last relation, (2.27), (2.17) and Lemma 2.3 lead to
\begin{equation}
12(5\alpha^{2} + \beta^{2})c + 64\alpha^{4} =
16\alpha^{2}(|W_{1}|^{2} + 3\beta^{2}) + 3 c^{2} .
\end{equation}
Because of (2.28): $f(\omega) = 64\omega^{2} + 60 c\omega  + 12 c\beta^{2}$, where $\omega = \alpha^{2}$, is positive for every $\omega, \beta$. This holds if and only if the
discriminant of $f(\omega)$ is negative for all $\beta, c$. But this is not true, hence we have a contradiction. Therefore $V'$ is empty and the real hypersurface $M$ consists only of $V_{0}$ and $V''$ i.e $\xi$ is principal
and $M$ is a Hopf hypersurface. \qed

\section{Proof of theorems}

From Lemma 2.4:
\begin{equation}
A\xi = \alpha\xi ,\hspace{50pt} \alpha = g(A\xi, \xi).
\end{equation}

We consider a $\phi - basis$ $\big\{ V_{i}, \phi V_{i}, \xi\big\} ,
(i = 1, 2, . . . n - 1)$. From (2.5) and (3.1) we obtain
\begin{equation}
lX = \frac{c}{4}\big[X - \eta(X)\xi \big] + \alpha AX -
\eta(X)\alpha^{2}\xi.
\end{equation}
(3.2) for $X = V_{i}$ implies
\begin{equation}
lV_{i} = \frac{c}{4}V_{i} + \alpha AV_{i}.
\end{equation}
Applying $\phi$ to (3.3) we obtain
\begin{equation}
\phi l V_{i} = \frac{c}{4}\phi V_{i} + \alpha \phi AV_{i}, \hspace{20pt} i = 1, . . . ,n - 1.
\end{equation}

 The relation (3.2) for $X = \phi V_{i}$
yields
\begin{equation}
l\phi V_{i} = \frac{c}{4}\phi V_{i} + \alpha A\phi V_{i}.
\end{equation}

Comparing (3.4) with (3.5), and by making use of the condition $\phi l = l
\phi$ we have
\begin{equation}
  (A\phi - \phi A)V_{i} = 0, \hspace{20pt} i = 1, . . . ,n - 1.
\end{equation}

On the other hand the action of $\phi$ on (3.5) yields
\begin{equation}
\phi (l\phi V_{i}) = \frac{c}{4}\phi^{2}V_{i} + \alpha \phi A\phi V_{i},
\end{equation}
which, by virtue of (1.1), is written in the form
\begin{equation}
(\phi l)\phi V_{i} = -\frac{c}{4}V_{i} + \alpha (\phi A)\phi V_{i},
\end{equation}
Moreover, the calculation of $(l\phi) \phi V_{i}$ by virtue of (1.1) and
(3.3) yields:
$$(l\phi) \phi V_{i} = l \phi^{2}V_{i} = - lV_{i} =
-\frac{c}{4}V_{i} - \alpha AV_{i} = -\frac{c}{4}V_{i} + \alpha A\phi^{2}
V_{i} = -\frac{c}{4}V_{i} + \alpha A\phi \phi V_{i} \Leftrightarrow$$
\begin{equation}
(l\phi) \phi V_{i} = -\frac{c}{4}V_{i} + \alpha (A\phi) \phi V_{i}
\end{equation}

Comparing (3.8) and (3.9), and by making use of the condition $\phi l = l
\phi$ we have
\begin{equation}
  (A\phi - \phi A)\phi V_{i} = 0
\end{equation}
 for every $i = 1, . . . ,n - 1$. But from (1.1) and (3.1) we also have
\begin{equation}
  (A\phi - \phi A)\xi = 0.
\end{equation}

So, (3.6), (3.10) and (3.11) imply that $A\phi = \phi A$. This
result and the Theorems (0.1) and (0.2) complete the proof of the main Theorems. \qed \\
\
\\

We must also notice that in class C (REMARK 1) we have $(\nabla_{\xi}l)V_{i} =
\mu\xi \Leftrightarrow \nabla_{\xi}lV_{i} - l\nabla_{\xi}V_{i}=
\mu\xi$, whose scalar product with $\xi$ yields $\mu = 0$. Also in
class D (REMARK 1) we have $(\nabla_{\xi}l)\xi = \mu\xi \Rightarrow \mu =
0$. So we have:

\begin{corollary}
Let M be a real hypersurface of a complex space form $M_{n}(c)$ $(c
\neq 0)$, satisfying $\phi l = l\phi$ and $\nabla_{\xi}l = \mu\xi$
on $ker(\eta)$ or on span$\{\xi\}$. Then the function $\mu$ must be
identically zero on M.
\end{corollary}
\bibliographystyle{amsplain}

\end{document}